# TRANSITION OPERATORS OF DIFFUSIONS REDUCE ZERO–CROSSING

STEVEN N. EVANS AND RUTH J. WILLIAMS

ABSTRACT. If $u(t, x)$ is a solution of a one–dimensional, parabolic, second–order, linear partial differential equation (PDE), then it is known that, under suitable conditions, the number of zero–crossings of the function $u(t, \cdot)$ decreases (that is, does not increase) as time $t$ increases. Such theorems have applications to the study of blow–up of solutions of semilinear PDE, time dependent Sturm Liouville theory, curve shrinking problems and control theory. We generalise the PDE results by showing that the transition operator of a (possibly time–inhomogeneous) one–dimensional diffusion reduces the number of zero–crossings of a function or even, suitably interpreted, a signed measure. Our proof is completely probabilistic and depends in a transparent manner on little more than the sample–path continuity of diffusion processes.

## 1. INTRODUCTION

The *number of zero–crossings* of a continuous function $f$ defined on some subinterval $I$ of $\mathbb{R}$ is the supremum of the set of $n \in \mathbb{N} := \{1, 2, \ldots\}$ for which there exist $x_1 < x_2 < \ldots < x_{n+1}$ in $I$ such that $f(x_i)f(x_{i+1}) < 0$ for $1 \leq i \leq n$ (with the convention that the supremum of the empty set is 0). An operator $K$ acting on some class of continuous functions on $I$ is said to be *variation diminishing* if the number of zero–crossings of $Kg$ is at most the number of zero–crossings of $g$ for all $g$ in the domain of $K$.

Consider a continuous solution of the one–dimensional, parabolic, second–order, linear initial value problem

$$(1.1) \quad \begin{aligned} \frac{\partial}{\partial s}u(s,x) &= \left(\frac{1}{2}a(s,x)\frac{\partial^2}{\partial x^2} + b(s,x)\frac{\partial}{\partial x} + c(s,x)\right)u(s,x), \ x \in I^o, \ s > 0, \\ u(0,x) &= f(x), \ x \in I, \end{aligned}$$

for some choice of boundary conditions (here $a(s,x) \geq 0$ and $I^o$ denotes the interior of $I$). The content of several results in the literature is that, under appropriate hypotheses on the coefficients and boundary conditions, the number of zero–crossings of $u(t, \cdot)$ is no more than the number of zero–crossings of $u(s, \cdot)$ when $0 \leq s < t$. Equivalently, by the flow property of solutions of (1.1), for each $t \geq 0$ the operator that takes the initial data $f$ to $u(t, \cdot)$ is variation diminishing. These results have

*Key words and phrases.* zero-crossing, variation diminishing, time-inhomogeneous diffusion, measure-valued process, martingale problem.

Evans' research was supported in part by NSF grant DMS-9703845. Williams' was supported in part by NSF grant DMS-9703891. Research at MSRI is supported in part by NSF grant DMS-9701755.

numerous applications that range from the study of blow–up of solutions of semilinear partial differential equations through time dependent Sturm Liouville theory and curve shrinking problems to control theory (see [21, 2, 3] for references). Many applications of the general notion of variation diminishing and its connections with total positivity can be found in [14]. For applications to statistics see [6].

Parabolic partial differential equations (PDE) and diffusion processes are intimately related (see [4, 8, 12, 22] for indications of some of the many connections), and our main goal in this paper is to show how results on the reduction of zero–crossings can follow rather simply from probabilistic considerations using little more than the fact that diffusion processes have continuous sample paths. The operators we consider are, in some ways, more general than those that take $f$ to $u(t, \cdot)$ via (1.1). In full generality, these operators act on signed measures rather than functions and, with an appropriate definition, the number of zero–crossings of a signed measure is diminished by their action. However, our results do not completely subsume those in the literature. For example, we always work with finite signed measures (and hence with integrable functions), whereas the result in [2] is stated just under a growth condition at infinity on the solution of (1.1). This finiteness restriction can be weakened somewhat under additional assumptions, but we do not pursue the matter here. Moreover, under suitable non-degeneracy hypotheses it is possible to obtain more precise information about the evolution of the zero set, such as the fact that double zeroes disappear immediately (see [2]). These finer results will typically fail to hold in the general setting we are considering.

Results on variation diminishing properties of operators that are associated with diffusions but do not necessarily come from PDE (and the connection with total positivity properties of such operators) can be found in [15, 16, 17, 13] and Ch. 5 of [14] (see also [20]). Our results differ from these in that we deal with the time–inhomogeneous case (which is necessary to subsume (1.1)), and we don't require that transition densities exist. Moreover, our proof is entirely probabilistic and, in particular, doesn't use the analytic technique of total positivity.

In order to describe our results, we need some preliminaries. We begin with a definition for the number of zero–crossings of a Radon signed measure on $I$ that coincides with that for continuous functions if we identify a continuous function $f$ with the measure that has Radon–Nikodym derivative $f$ against Lebesgue measure.

**Definition 1.1.** Given a Radon signed measure $\mu$ on $I$ with Hahn–Jordan decomposition $\mu = \mu^+ - \mu^-$, define Radon non-negative measures $\mu_n^+$ and $\mu_n^-$ on the Cartesian product $I^n$ for each $n \in \mathbb{N}$ as follows:

$$\mu_1^+ := \mu^+, \ \mu_2^+ := \mu^+ \otimes \mu^-, \ \mu_3^+ := \mu^+ \otimes \mu^- \otimes \mu^+, \ \ldots$$
$$\mu_1^- := \mu^-, \ \mu_2^- := \mu^- \otimes \mu^+, \ \mu_3^- := \mu^- \otimes \mu^+ \otimes \mu^-, \ \ldots$$

Write $H_n := \{(x_1, \ldots, x_n) \in I^n : x_1 < x_2 < \ldots < x_n\}$ and define $\Sigma(\mu)$, the *number of zero–crossings of $\mu$* by

$$\Sigma(\mu) := \sup\{n : \mu_n^+(H_n) > 0 \text{ or } \mu_n^-(H_n) > 0\} - 1.$$

We now describe our generalisations for the operators that take the initial data $f$ to the solution $u(t, \cdot)$ at time $t \geq 0$ in (1.1). Put $E := \mathbb{R}_+ \times I$, where $\mathbb{R}_+ = [0, \infty[$. Fix a strongly continuous, positive, contraction semigroup $(P_t)_{t \geq 0}$ of linear operators on $C_0(E)$ ( $\equiv$ continuous functions on $E$ that "vanish at $\infty$"). Such a semigroup is usually called Feller or Feller–Dynkin in the probability literature, and we refer the reader to Chapter III of [19] or Chapters 1 and 4 of [9] for the standard Hille–Yosida theory and facts about the Markov process associated with such semigroups that we use below. Denote by $A$ the generator of $(P_t)_{t \geq 0}$, and write $\mathcal{D}(A)$ for the domain of $A$. With the usual slight abuse of notation, let $P_t(w, \cdot)$, $w \in E$, denote the subprobability measure on $E$ such that $(P_t f)(w) = \int_E P_t(w, dw') f(w')$, $f \in C_0(E)$. Assume that $(P_t)_{t \geq 0}$ is a *space-time* semigroup; that is, the measure $P_t((s, x), \cdot)$ is concentrated on $\{s + t\} \times I$.

Let $\partial$ be the point at infinity in the one–point compactification of $E$, denoted by $E^\partial$. Define a family $P_t^\partial(w, \cdot)$, $t \geq 0$, $w \in E^\partial$, of probability measures on $E^\partial$ by

$$P_t^\partial(w, B) := P_t(w, B \cap E) + (1 - P_t(w, E))\delta_\partial(B), \ w \in E,$$
$$P_t^\partial(\partial, B) := \delta_\partial(B),$$

where $\delta_\partial$ denotes the point mass at $\partial$. With the same abuse of notation used above, the family of operators $(P_t^\partial f)(w) = \int_{E^\partial} P_t^\partial(w, dw') f(w')$, $f \in C(E^\partial)$, defines a strongly continuous, positive, contraction semigroup on $C(E^\partial)$, and there is a strong Markov process $(W, (\mathbb{P}^w)_{w \in E^\partial})$ with state–space $E^\partial$ and transition semigroup $(P_t^\partial)_{t \geq 0}$. The state $\partial$ has the property that $W_t = \partial$ for all $t \geq \zeta := \inf\{s \geq 0 : W_s = \partial\}$; and $\partial$ and $\zeta$ are called the *cemetery* and *life–time*, respectively. Denote by $A^\partial$ the generator of $(P_t^\partial)_{t \geq 0}$, and write $\mathcal{D}(A^\partial)$ for the domain of $A^\partial$.

Adjoin another fictitious state $\dagger$ to $I$ to form $I^\dagger$. If $I$ is compact, then take $\dagger$ to be an isolated point. Otherwise, $\dagger$ is the point at infinity in the one–point compactification of $I$. It will simplify notation to interpret $(s, \dagger)$ as $\partial$ for all $s \geq 0$. With this interpretation, given a function $f : E^\partial \to \mathbb{R}$ define $f_s : I^\dagger \to \mathbb{R}$, $s \geq 0$, by $f_s = f(s, \cdot)$.

Write $W_t =: (T_t, X_t) \in \mathbb{R}_+ \times I$ for $0 \leq t < \zeta$. Put $X_t := \dagger$ when $t \geq \zeta$. Suppose that $W$ is a *diffusion process*; that is, the sample paths of $W$ are continuous on $[0, \zeta[$. There appears to be no known analytic necessary and sufficient condition on $(P_t)_{t \geq 0}$ for $W$ to be a diffusion; however, a sufficient condition is that for all compact sets $K \subset E$ and all open neighbourhoods $U$ such that $K \subset U$, $\lim_{t \downarrow 0} \sup_{w \in K} t^{-1} P_t(w, E \setminus U) = 0$ (see Proposition I.9.10 of [5]). The assumption that $(P_t)_{t \geq 0}$ is a space-time semigroup is equivalent to the statement that $T_t = s + t$ for all $0 \leq t < \zeta$ almost surely under $\mathbb{P}^{(s,x)}$, $(s, x) \in E$.

Define a family $(Q_t)_{t \geq 0}$ of subprobability kernels on $I$ by $Q_t(x, B) := P_t((0, x), \mathbb{R}_+ \times B) = P_t((0, x), \{t\} \times B)$. As usual, if $\mu$ is a finite signed measure on $I$, define another finite signed measure $\mu Q_t$ on $I$ by $(\mu Q_t)(B) := \int_I \mu(dx) Q_t(x, B)$. Define a family $(Q_t^\partial)_{t \geq 0}$ of probability kernels on $I^\dagger$ from $(P_t^\partial)_{t \geq 0}$ in an analogous manner.

Our main result is the following.

**Theorem 1.2.** *Under the above assumptions, $\Sigma(\mu Q_s) \geq \Sigma(\mu Q_t)$ for all finite signed measures $\mu$ on $I$ and all $0 \leq s < t < \infty$.*

The most important class of examples of this set–up is when $I = \mathbb{R}$ and the generator $A$ extends the second order differential operator

$$\frac{\partial}{\partial s} + \frac{1}{2}\alpha(s,x)\frac{\partial^2}{\partial x^2} + \beta(s,x)\frac{\partial}{\partial x} + \gamma(s,x),$$

where $\alpha \geq 0$ and $\gamma \leq 0$. Under appropriate conditions, for example, if $\alpha$, $\beta$, $\gamma$, $\frac{\partial \alpha}{\partial x}$, $\frac{\partial \beta}{\partial x}$, $\frac{\partial^2 \alpha}{\partial x^2}$ are all bounded and uniformly Hölder continuous, $\alpha \geq c$ for some constant $c > 0$, and $\mu(A) := \int_A f(x)dx$ for a bounded, integrable, continuous function $f$ (cf. Section 5.7B of [12], Chapter 1 of [10] and Sections 6.4–6.5 of [11]), we have that for each $s > 0$, $\mu Q_s$ has a Radon–Nikodym derivative $u(s,\cdot)$ with respect to Lebesgue measure that solves

(1.2)
$$\frac{\partial}{\partial s}u(s,x) = \frac{1}{2}\frac{\partial^2}{\partial x^2}\left(\alpha(s,x)u(s,x)\right) - \frac{\partial}{\partial x}\left(\beta(s,x)u(s,x)\right) + \gamma(s,x)u(s,x), \quad s > 0,$$
$$\lim_{s \downarrow 0} u(s,x) = f(x).$$

If

$$\alpha(s,x) = a(s,x)$$
$$\beta(s,x) = \frac{\partial}{\partial x}a(s,x) - b(s,x)$$

and

$$\gamma(s,x) = \frac{1}{2}\frac{\partial^2}{\partial x^2}a(s,x) - \frac{\partial}{\partial x}b(s,x) + c(s,x),$$

then on setting $u(0,\cdot) = f(\cdot)$, we have that $u$ is the continuous solution on $\mathbb{R}_+ \times \mathbb{R}$ of (1.1). Of course, it may not be the case that this choice of $\gamma$ is non-positive, but this difficulty can be resolved when the choice of $\gamma$ is bounded by noting that for any constant $\kappa$ the measures $e^{\kappa s}\mu Q_s$ have Radon–Nikodym derivatives with respect to Lebesgue measure that solve (1.2) with $\gamma$ replaced by $\gamma + \kappa$.

The technique we use for proving Theorem 1.2 is to approximate the (deterministic) signed measure-valued function $(\mu Q_t)_{t \geq 0}$ by a discrete signed measure-valued stochastic process for which the analogue of Theorem 1.2 is obvious. A sketch of the idea of the proof is as follows. Assume without loss of generality that $|\mu| := \mu^+ + \mu^-$ is a probability measure and for $n \in \mathbb{N}$ let $(\mathcal{E}_1, \mathcal{X}_1), \ldots, (\mathcal{E}_n, \mathcal{X}_n)$ be i.i.d. $\{+1,-1\} \times I$-valued random variables with

(1.3)    $\mathbb{P}\{\mathcal{E}_i = +1, \mathcal{X}_i \in B\} = \mu^+(B)$ and $\mathbb{P}\{\mathcal{E}_i = -1, \mathcal{X}_i \in B\} = \mu^-(B).$

Conditional on $\mathcal{X}_i = x_i$ for $1 \leq i \leq n$ build $n$ independent processes $X^1, \ldots, X^n$, with $X^i$ having the law of $X$ under $\mathbb{P}^{(0,x_i)}$. Think of $X^i$ as describing the motion of a particle with sign $\mathcal{E}_i$. When particles of opposite sign meet, let them annihilate in pairs. Write $C_t$ for the set of indices of the particles that at time $t$ have not been annihilated or sent to †. It is almost immediate that $\Sigma(\sum_{i \in C_t} \mathcal{E}_i \delta_{X_t^i})$ is non-increasing in $t$. Using the law of large numbers and the martingale problem characterisation of the law of $X$, we show that $n^{-1}\sum_{i \in C_t} \mathcal{E}_i \delta_{X_t^i}$ converges in probability to $\mu Q_t$ as $n \to \infty$. The fact that the zero-crossings of $n^{-1} \sum_{i \in C_t} \mathcal{E}_i \delta_{X_t^i}$ do not increase with time is inherited by $\mu Q_t$ and this implies the theorem.

The plan of the rest of the paper is as follows. In Section 2 we collect together some facts about signed measures and give two useful equivalent definitions of the number of zero–crossings. We introduce the approximating process of signed particles in Section 3 and develop some of its properties. We complete the proof of Theorem 1.2 in Section 4. In Section 5 we exhibit a simple counterexample to demonstrate that the assumption of sample–path continuity of $W$ up to the life–time cannot be removed.

**Notation 1.3.** We will often use the functional notation $\nu(g) := \int g d\nu$ for integrals. Similarly, we will use the notation $\mathbb{P}[V]$ to denote the expectation of a random variable $V$ with respect to a probability measure $\mathbb{P}$.

## 2. Signed measures and sign sequences

The following simple result almost certainly appears in the literature, but we have been unable to find a reference.

**Lemma 2.1.** *Let $\mu$ be a Radon signed measure on $\mathbb{R}$. Suppose that $\mu^+([a,b]) > 0$ for $-\infty < a \leq b < \infty$. Then for every $\epsilon > 0$ there exists a continuous, non-negative function $g$ with support contained in $[a - \epsilon, b + \epsilon]$ such that $\mu(g) > 0$.*

*Proof.* By the Hahn–Jordan decomposition, there exists a Borel set $B \subseteq [a, b]$ such that $\mu^+(B) > 0$ and $\mu^-(B) = 0$. By inner regularity, there exists a compact set $K \subseteq B$ such that $\mu^+(K) > 0$ (and, of course, $\mu^-(K) = 0$). By outer regularity, there exists an open set $U$ with $K \subseteq U \subseteq [a - \epsilon, b + \epsilon]$ such that $\mu^-(U) < \mu^+(K)$. Now take $g$ to have values in $[0, 1]$, to be $1$ on $K$ and to be $0$ on $\mathbb{R} \setminus U$. Then $\mu(g) = \mu^+(g) - \mu^-(g) \geq \mu^+(K) - \mu^-(U) > 0$, as required. $\square$

In the following result and later, it will be convenient to identify a Radon signed measure $\mu$ on $I$ with the Radon signed measure $\tilde{\mu}$ on $\mathbb{R}$ such that $\tilde{\mu}(B) = \mu(B)$ for $B \subseteq I$ and $\tilde{\mu}(\mathbb{R} \setminus I) = 0$.

**Lemma 2.2.** *For a non-zero Radon signed measure $\mu$ on $I$, the number of zero–crossings $\Sigma(\mu)$ is the supremum of the set of $n \in \mathbb{N}$ for which there are continuous, non-negative, compactly supported functions $g_1, \ldots, g_{n+1}$ defined on $\mathbb{R}$ such that:*

*(a) $g_i(x) > 0$ and $g_j(y) > 0$ for $i < j$ only if $x < y$,*
*(b) $\mu(g_k)\mu(g_{k+1}) < 0$ for $1 \leq k \leq n$,*
*(where we adopt the convention that the supremum of the empty set is $0$).*

*Proof.* Let $\Sigma'(\mu)$ denote the quantity defined in the statement. If $\Sigma(\mu) = 0$, then it is certainly the case that $\Sigma'(\mu) \geq \Sigma(\mu)$. Suppose that $\Sigma(\mu) \geq n$ for some $n \in \mathbb{N}$. Assume without loss of generality that, in the notation used in the definition of $\Sigma(\mu)$, $\mu_{n+1}^+(H_{n+1}) > 0$. By inner regularity, there exists $-\infty < a_1 \leq b_1 < a_2 \leq b_2 < \ldots < a_{n+1} \leq b_{n+1} < \infty$ such that $\mu_{n+1}^+([a_1, b_1] \times \cdots \times [a_{n+1}, b_{n+1}]) > 0$ and hence $\mu^+([a_1, b_1]) > 0, \mu^-([a_2, b_2]) > 0, \ldots$. It is clear how to use Lemma 2.1 and its counterpart for $\mu^-$ to construct continuous, non-negative, compactly supported functions $g_1, \ldots, g_{n+1}$ satisfying (a) and (b), and hence $\Sigma'(\mu) \geq n$. Thus $\Sigma'(\mu) \geq \Sigma(\mu)$.

For the reverse inequality, if $\Sigma'(\mu) = 0$, then it is again certainly the case that $\Sigma(\mu) \geq \Sigma'(\mu)$. Suppose that $\Sigma'(\mu) \geq n$ for some $n \in \mathbb{N}$. Let $g_1, \ldots, g_{n+1}$ be as in the statement of the lemma. Assume without loss of generality, that $\mu(g_1) > 0$, so that $\mu(g_2) < 0$, $\mu(g_3) > 0$, and so on. It follows that there are Borel sets $A_1, \ldots, A_{n+1}$ with $x \in A_i$ and $y \in A_j$ for $i < j$ only if $x < y$ and such that $\mu^+(A_1) > 0$, $\mu^-(A_2) > 0$, and so on. Clearly, $\mu_{n+1}^+(H_{n+1}) \geq \mu_{n+1}^+(A_1 \times \cdots \times A_{n+1}) > 0$, and hence $\Sigma(\mu) \geq n$. Thus $\Sigma(\mu) \geq \Sigma'(\mu)$. □

A *sign sequence* is either the empty sequence $\emptyset$ or a finite length sequence $(\epsilon_1, \ldots, \epsilon_n)$, where $\epsilon_i \in \{+1, -1\}$. Define a partial order on sign sequences by saying that $a \preccurlyeq b$ if $a$ is a subsequence of $b$ (where the empty sequence is a subsequence of every sequence). For a sign sequence $a = (\epsilon_1, \ldots, \epsilon_n)$, let

$$\sigma(a) := \max\{k \in \mathbb{N} : \exists i_1 < \ldots < i_{k+1}, \epsilon_{i_j} \neq \epsilon_{i_{j+1}}, 1 \leq j \leq k\}$$

(with the convention that the maximum of the empty set is 0) and put $\sigma(\emptyset) := -1$. Observe that if $a \preccurlyeq b$, then $\sigma(a) \leq \sigma(b)$. A *substring* of a sign sequence $a = (\epsilon_1, \ldots, \epsilon_n)$ is a subsequence consisting of consecutive terms, that is, one of the form $(\epsilon_u, \epsilon_{u+1}, \ldots, \epsilon_{v-1}, \epsilon_v)$.

Let $\mathcal{M}$ and $\mathcal{N}$ denote, respectively, the space of finite signed measures on $I$ and the space of finite signed integer–valued measures on $I$. Equip $\mathcal{M}$ and $\mathcal{N}$ with the weak topology. Let $\mathcal{M}^+$ and $\mathcal{N}^+$ denote the subsets of non-negative measures in $\mathcal{M}$ and $\mathcal{N}$, respectively. Write $\overline{\mathcal{M}}$ for the space of finite signed measures on $I^\dagger$, again equipped with the weak topology. Define $\overline{\mathcal{N}}, \overline{\mathcal{M}}^+$ and $\overline{\mathcal{N}}^+$ in the obvious way. Extend the definition of $\Sigma$ to $\overline{\mathcal{M}}$ by setting $\Sigma(\bar{\mu}) = \Sigma(\bar{\mu}_{|I})$. When there is no danger of confusion, identify $\mu \in \mathcal{M}$ with the measure $\mu(\cdot \cap I) \in \overline{\mathcal{M}}$.

Given $\nu \in \mathcal{N} \backslash \{0\}$, represent $\nu$ as

$$\nu = \sum_{i=1}^{n} \epsilon_i \delta_{x_i},$$

where $\epsilon_i \in \{+1, -1\}$, $x_1 \leq \ldots \leq x_n$ and $x_i < x_{i+1}$ if $\epsilon_i \neq \epsilon_{i+1}$, and write $\mathcal{S}(\nu) := (\epsilon_1, \ldots, \epsilon_n)$. Put $\mathcal{S}(0) := \emptyset$. Extend the definition of $\mathcal{S}$ to $\overline{\mathcal{N}}$ by setting $\mathcal{S}(\bar{\nu}) := \mathcal{S}(\bar{\nu}_{|I})$. It is clear that for all $\bar{\nu} \in \overline{\mathcal{N}}$,

(2.1) $$\Sigma(\bar{\nu}) = \sigma \circ \mathcal{S}(\bar{\nu}).$$

## 3. A SIGNED MEASURE–VALUED PROCESS

Fix, for the moment, $\nu \in \mathcal{N}$. In the next few paragraphs we define a càdlàg $\overline{\mathcal{N}} \times \overline{\mathcal{N}}^+$-valued process $((Y_t, Z_t))_{t \geq 0}$ with $(Y_0, Z_0) = (\nu, |\nu|)$. The process $Z$ can be thought of as a system of unsigned particles that move as independent copies of $X$. In constructing $Y$ the same particles are given signs, and when a pair of $+1$ and $-1$ particles collide at a point other than the cemetery $\dagger$ they annihilate each other. If $X$ is a "non-degenerate" process such as Brownian motion, then it is almost surely the case that particles only collide in pairs. However, because of the generality in which we are working we need to allow for the possibility of collisions involving more particles.

Write $\nu = \sum_{i=1}^{n^+} \delta_{x^{+,i}} - \sum_{j=1}^{n^-} \delta_{x^{-,j}}$, where it may be that $x^{+,i} = x^{+,i'}$ for $i \ne i'$ and $x^{-,j} = x^{-,j'}$ for $j \ne j'$ but $x^{+,k} \ne x^{-,\ell}$ for all $k, \ell$. Let $(X_t^{+,i})_{t \ge 0}$, $1 \le i \le n^+$, and $(X_t^{-,j})_{t \ge 0}$, $1 \le j \le n^-$, be independent $I^\dagger$-valued processes with $X^{+,i}$ (resp. $X^{-,j}$) having the law of $X$ under $\mathbb{P}^{(0,x^{+,i})}$ (resp. $\mathbb{P}^{(0,x^{-,j})}$).

Define stopping times $0 =: S_0 \le S_1 \le \ldots$ and random sets of indices $\{1, \ldots, n^+\} =: J_0^+ \supseteq J_1^+ \supseteq \ldots$ and $\{1, \ldots, n^-\} =: J_0^- \supseteq J_1^- \supseteq \ldots$ inductively as follows. Put

$$S_{m+1} := \inf\{t > S_m : X_t^{+,k} = X_t^{-,\ell} \ne \dagger \text{ for some } k \in J_m^+ \text{ and } \ell \in J_m^-\}.$$

If $S_{m+1} < \infty$, let $x_1 < x_2 < \ldots < x_r \in I$ be a listing of the random set of points $x$ such that there is at least one $k \in J_m^+$ and one $\ell \in J_m^-$ for which $X_{S_{m+1}}^{+,k} = X_{S_{m+1}}^{-,\ell} = x$ (for simplicity, we suppress the dependence on $m$ of this and other notation involved in the definition of $J_{m+1}^+$ and $J_{m+1}^-$). Put

$$K_q^+ := \{k \in J_m^+ : X_{S_{m+1}}^{+,k} = x_q\}$$

and

$$K_q^- := \{\ell \in J_m^- : X_{S_{m+1}}^{-,\ell} = x_q\}.$$

Set $G_q := \#(K_q^+) \wedge \#(K_q^-)$. Write $L_q^+$ (resp. $L_q^-$) for the $G_q$ largest elements of $K_q^+$ (resp. $K_q^-$). Put

$$J_{m+1}^+ := J_m^+ \setminus \bigcup_{q=1}^r L_q^+$$

and

$$J_{m+1}^- := J_m^- \setminus \bigcup_{q=1}^r L_q^-$$

If $S_{m+1} = \infty$, let $J_{m+1}^+ = J_m^+$ and $J_{m+1}^- = J_m^-$.

The times $S_1, S_2, \ldots$ that are finite are the distinct times at which collisions between $+1$ and $-1$ particles occur. If $S_m < \infty$, then $J_m^+$ (resp. $J_m^-$) is the set of indices of $+1$ (resp. $-1$) particles still alive at time $S_m$. Each time a collision between $+1$ and $-1$ particles occurs, the number of $+1$ particles and the number of $-1$ particles are both reduced by at least one. That is, if $S_m < \infty$, then $\#(J_m^+) < \#(J_{m-1}^+)$ and $\#(J_m^-) < \#(J_{m-1}^-)$. Therefore, $S_M = \infty$ for $M$ sufficiently large ($M > n^+ \wedge n^-$ works).

Now define

(3.1) $$Y_t := \sum_{i \in J_m^+} \delta_{X_t^{+,i}} - \sum_{j \in J_m^-} \delta_{X_t^{-,j}}, \; S_m \le t < S_{m+1},$$

and

$$Z_t := \sum_{i=1}^{n^+} \delta_{X_t^{+,i}} + \sum_{j=1}^{n^-} \delta_{X_t^{-,j}}. \tag{3.2}$$

Put $\tau^{+,i} := \inf\{S_m : i \notin J_m^+\}$, $1 \leq i \leq n^+$, and $\tau^{-,j} := \inf\{S_m : j \notin J_m^-\}$, $1 \leq j \leq n^-$. Set $\hat{X}_t^{+,i} := X_{t \wedge \tau^{+,i}}^{+,i}$ and $\hat{X}_t^{-,j} := X_{t \wedge \tau^{-,j}}^{-,j}$. Then,

$$Y_t = \sum_{i=1}^{n^+} \delta_{\hat{X}_t^{+,i}} - \sum_{j=1}^{n^-} \delta_{\hat{X}_t^{-,j}}, \tag{3.3}$$

since, for the construction of the signed measure $Y_t$, freezing pairs of particles of opposite signs where they collide has the same effect as annihilating them.

**Lemma 3.1.** *(i) For all $t \geq 0$, $|Y_t| \leq Z_t$.*
*(ii) For all $t \geq 0$, $Z_t(I^\dagger) = |\nu|(I)$.*
*(iii) For any $f \in \mathcal{D}(A^\partial)$ the processes*

$$M_t^{Y,f} := Y_t(f_t) - Y_0(f_0) - \int_0^t Y_s((A^\partial f)_s)\, ds$$

*and*

$$M_t^{Z,f} := Z_t(f_t) - Z_0(f_0) - \int_0^t Z_s((A^\partial f)_s)\, ds$$

*are càdlàg martingales that are bounded on finite intervals and such that*

$$[M^{Y,f}]_t - [M^{Y,f}]_s \leq [M^{Z,f}]_t - [M^{Z,f}]_s,\ 0 \leq s \leq t,$$

*where $[\cdot]$ denotes quadratic variation.*
*(iv) The process $(\Sigma(Y_t))_{t \geq 0}$ has non-increasing sample paths almost surely.*

*Proof.* Parts (i) and (ii) are obvious. Consider part (iii). By Proposition 4.1.7 of [9], the processes

$$M_t^{+,i,f} := f(t, X_t^{+,i}) - f(0, x^{+,i}) - \int_0^t A^\partial f(s, X_s^{+,i})\, ds$$

and

$$M_t^{-,j,f} := f(t, X_t^{-,j}) - f(0, x^{-,j}) - \int_0^t A^\partial f(s, X_s^{-,j})\, ds$$

are càdlàg martingales. By (3.2) and (3.3),

$$M_t^{Z,f} = \sum_{i=1}^{n^+} M_t^{+,i,f} + \sum_{j=1}^{n^-} M_t^{-,j,f},$$

and
$$M_t^{Y,f} = \sum_{i=1}^{n^+} M_{t\wedge\tau^{+,i}}^{+,i,f} - \sum_{j=1}^{n^-} M_{t\wedge\tau^{-,j}}^{-,j,f},$$

where in the last line we have taken into account cancellation of like terms from the two sums, which occurs whenever $\tau^{+,i} = \tau^{-,j} < t$ and $X_{\tau^{+,i}}^{+,i} = X_{\tau^{-,j}}^{-,j}$ for some $i$ and $j$. It is clear from the above that both $M^{Z,f}$ and $M^{Y,f}$ are locally bounded martingales.

Consider the claim regarding quadratic variations. Note that

(3.4) $\quad [M_{\cdot\wedge\tau^{+,i}}^{+,i,f}]_t = [M^{+,i,f}]_{t\wedge\tau^{+,i}}$ and $[M_{\cdot\wedge\tau^{-,j}}^{-,j,f}]_t = [M^{-,j,f}]_{t\wedge\tau^{-,j}}.$

By construction, for each $t \geq 0$ the collection $\{M^{+,i,f} : 1 \leq i \leq n^+\} \cup \{M^{-,j,f} : 1 \leq j \leq n^-\}$ is conditionally independent given $\{M_s^{+,i,f} : 1 \leq i \leq n^+, 0 \leq s \leq t\} \cup \{M_s^{-,j,f} : 1 \leq j \leq n^-, 0 \leq s \leq t\}$. Consequently,

(3.5)
$$[M^{+,i,f}, M^{+,i',f}] \equiv 0, \quad i \neq i',$$
$$[M^{-,j,f}, M^{-,j',f}] \equiv 0, \quad j \neq j',$$
$$[M^{+,i,f}, M^{-,j,f}] \equiv 0.$$

Thus,

(3.6)
$$[M_{\cdot\wedge\tau^{+,i}}^{+,i,f}, M_{\cdot\wedge\tau^{+,i'}}^{+,i',f}] \equiv 0, \quad i \neq i',$$
$$[M_{\cdot\wedge\tau^{-,j}}^{-,j,f}, M_{\cdot\wedge\tau^{-,j'}}^{-,j',f}] \equiv 0, \quad j \neq j',$$
$$[M_{\cdot\wedge\tau^{+,i}}^{+,i,f}, M_{\cdot\wedge\tau^{-,j}}^{-,j,f}] \equiv 0,$$

(by, for example, combining parts (i) and (vii) of Theorem VI.37.9 [18]).

Therefore,

$$[M^{Y,f}]_t - [M^{Y,f}]_s = \sum_{i=1}^{n^+}\left\{[M_{\cdot\wedge\tau^{+,i}}^{+,i,f}]_t - [M_{\cdot\wedge\tau^{+,i}}^{+,i,f}]_s\right\} + \sum_{j=1}^{n^-}\left\{[M_{\cdot\wedge\tau^{-,j}}^{-,j,f}]_t - [M_{\cdot\wedge\tau^{-,j}}^{-,j,f}]_s\right\}$$

$$= \sum_{i=1}^{n^+}\left\{[M^{+,i,f}]_{t\wedge\tau^{+,i}} - [M^{+,i,f}]_{s\wedge\tau^{+,i}}\right\} + \sum_{j=1}^{n^-}\left\{[M^{-,j,f}]_{t\wedge\tau^{-,j}} - [M^{-,j,f}]_{s\wedge\tau^{-,j}}\right\}$$

$$\leq \sum_{i=1}^{n^+}\left\{[M^{+,i,f}]_t - [M^{+,i,f}]_s\right\} + \sum_{j=1}^{n^-}\left\{[M^{-,j,f}]_t - [M^{-,j,f}]_s\right\}$$

$$= [M^{Z,f}]_t - [M^{Z,f}]_s,$$

where the successive equalities follow by (3.6), (3.4) and (3.5), respectively.

Turning to part (iv), it follows from the construction and the continuity of sample paths of $X$ up to the life–time $\zeta$ that

(3.7) $\quad\quad\quad\quad\quad\quad\quad\quad \mathcal{S}(Y_t) \preccurlyeq \mathcal{S}(Y_s),\ 0 \leq s < t < \infty.$

More specifically, if $S_m \leq s < t < S_{m+1}$ for some $m$, then $\mathcal{S}(Y_t)$ is obtained from $\mathcal{S}(Y_s)$ by the removal of elements corresponding to particles that are killed in the interval $]s,t]$. If $S_m \leq s < t = S_{m+1}$ for some $m$, then $\mathcal{S}(Y_t)$ is obtained from $\mathcal{S}(Y_s)$ by the removal of elements corresponding to killed particles plus the possibly repeated substitution (due to the mutual annihilation of oppositely signed pairs of particles) of a substring containing both $+1$ and $-1$ elements by either an empty substring (if the numbers of $+1$ and $-1$ particles colliding at a point are equal) or a substring with all elements the same (if the numbers colliding are not equal). Thus (3.7) holds when $S_m \leq s < t \leq S_{m+1}$ for some $m$, and the case of general $0 \leq s < t < \infty$ holds by transitivity. Therefore, by (2.1),

$$\Sigma(Y_t) = \sigma \circ \mathcal{S}(Y_t) \leq \sigma \circ \mathcal{S}(Y_s) = \Sigma(Y_s),\ 0 \leq s < t < \infty.$$

□

## 4. Proof of Theorem 1.2

Suppose without loss of generality that $|\mu| := \mu^+ + \mu^-$ is a probability measure. Define probability measures $\Pi^n$, $n \in \mathbb{N}$, on $\mathcal{N}$ by

$$\Pi^n(F) := \int_{(\{+1,-1\} \times I)^n} \left(\delta_{+1} \otimes \mu^+ + \delta_{-1} \otimes \mu^-\right)^{\otimes n}(d(\epsilon,x))\ F\left(\sum_{i=1}^n \epsilon_i \delta_{x_i}\right)$$

for $F$ a bounded Borel function on $\mathcal{N}$. That is, $\Pi^n$ is the law of the random signed measure $\sum_{i=1}^n \mathcal{E}_i \delta_{\mathcal{X}_i}$ described in (1.3). Observe that

(4.1) $$\Pi^n\left(\{\nu : \Sigma(\nu) > \Sigma(\mu)\}\right) = 0.$$

For each $\nu \in \mathcal{N}$, let $\mathbb{P}^\nu$ denote the law of the process $(Y,Z)$ constructed in Section 3. It is clear that on some probability space $(\Omega, \mathcal{F}, \mathbb{P})$ it is possible to build a sequence $((Y^n, Z^n))_{n=1}^\infty$ of $\overline{\mathcal{N}} \times \overline{\mathcal{N}}^+$-valued processes such that:

- the pair $(Y^n, Z^n)$ has law $\int \Pi^n(d\nu) \mathbb{P}^\nu$,
- the sequence of random measures $Y_0^1, Y_0^2 - Y_0^1, Y_0^3 - Y_0^2, \ldots$ is i.i.d.,
- the sequence of processes $Z^1, Z^2 - Z^1, Z^3 - Z^2, \ldots$ is i.i.d.

Put $\tilde{Y}^n := n^{-1} Y^n$ and $\tilde{Z}^n := n^{-1} Z^n$. The proof of Theorem 1.2 hinges on the following result.

**Lemma 4.1.** *The sequence $(\tilde{Y}^n)_{n=1}^\infty$ of càdlàg $\overline{\mathcal{M}}$–valued processes converges in probability in the Skorohod topology to the continuous deterministic function $(\mu Q_t^\partial)_{t \geq 0}$.*

*Proof.* We begin by establishing that the sequence $(\tilde{Y}^n)_{n=1}^\infty$ is tight and hence, by Prohorov's theorem, relatively compact in distribution. This will follow if we can prove that the sequence of finite signed measure-valued processes $((\delta_t \otimes \tilde{Y}_t^n)_{t \geq 0})_{n=1}^\infty$ is tight. As $\mathcal{D}(A^\partial)$ is a dense subspace of $C(E^\partial)$ (see, for example, Theorem 1.2.6 of [9]), an adaptation of the proof of Theorem 3.7.1 of [7] gives that it is enough to show for each $f \in \mathcal{D}(A^\partial)$ that the sequence of càdlàg real-valued processes $((\tilde{Y}_t^n(f_t))_{t \geq 0})_{n=1}^\infty$ is tight.

Fix $f \in \mathcal{D}(A^\partial)$. By the strong law of large numbers and the fact that $Y_0^1, Y_0^2 - Y_0^1, Y_0^3 - Y_0^2, \ldots$ is i.i.d., the uniformly bounded sequence $(\tilde{Y}_0^n(f_0))_{n=1}^\infty$ converges almost surely to $\mu(f_0)$.

From part (iii) of Lemma 3.1, the processes

$$M_t^{\tilde{Y}^n, f} := \tilde{Y}_t^n(f_t) - \tilde{Y}_0^n(f_0) - \int_0^t \tilde{Y}_s^n((A^\partial f)_s)\, ds$$

and

$$M_t^{\tilde{Z}^n, f} := \tilde{Z}_t^n(f_t) - \tilde{Z}_0^n(f_0) - \int_0^t \tilde{Z}_s^n((A^\partial f)_s)\, ds$$

are càdlàg martingales that are bounded on finite intervals and such that

$$[M^{\tilde{Y}^n, f}]_t - [M^{\tilde{Y}^n, f}]_s \leq [M^{\tilde{Z}^n, f}]_t - [M^{\tilde{Z}^n, f}]_s,\ 0 \leq s \leq t.$$

Consider a bounded stopping time $\tau$ and a constant $\eta > 0$. It follows from parts (i) and (ii) of Lemma 3.1 that $\sup_n \sup_t |\tilde{Y}_t^n|(I^\dagger) \leq 1$ and hence

$$\left| \int_\tau^{\tau+\eta} \tilde{Y}_s^n((A^\partial f)_s)\, ds \right| \leq \|A^\partial f\|_\infty \eta.$$

It follows from part (iii) of Lemma 3.1 and the fact that $Z^1, Z^2 - Z^1, \ldots$ are i.i.d. that

(4.2)
$$\begin{aligned}
\mathbb{P}\left[\left(M_{\tau+\eta}^{\tilde{Y}^n,f} - M_\tau^{\tilde{Y}^n,f}\right)^2\right] &= \mathbb{P}\left[[M^{\tilde{Y}^n,f}]_{\tau+\eta} - [M^{\tilde{Y}^n,f}]_\tau\right] \\
&\leq \mathbb{P}\left[[M^{\tilde{Z}^n,f}]_{\tau+\eta} - [M^{\tilde{Z}^n,f}]_\tau\right] \\
&= n^{-1} \mathbb{P}\left[[M^{\tilde{Z}^1,f}]_{\tau+\eta} - [M^{\tilde{Z}^1,f}]_\tau\right].
\end{aligned}$$

A standard application of Aldous's criterion for tightness (see [1]) completes the proof of tightness.

It follows from (4.2) and Doob's $L^2$-maximal inequality that

$$\lim_{n\to\infty} \mathbb{P}\left[\sup_{0\leq s\leq t} \left(M_s^{\tilde{Y}^n,f}\right)^2\right] = 0,\ \text{for all } t \geq 0.$$

Hence, if $\tilde{Y}^\infty$ is any subsequential limit in distribution of $(\tilde{Y}^n)_{n=1}^\infty$, then

$$\tilde{Y}_t^\infty(f_t) = \mu(f_0) + \int_0^t \tilde{Y}_s^\infty((A^\partial f)_s)\, ds \text{ for all } t \geq 0,$$

using the right continuity of $\tilde{Y}_\cdot^\infty(f_\cdot)$ and the fact that convergence in distribution in the Skorohod topology implies convergence of finite dimensional distributions for all but a countable set of times (see Theorem 7.8 of [9]).

Note also that $\sup_t |\tilde{Y}_t^\infty|(I^\dagger) \leq 1$. As $(\lambda - A^\partial)(\mathcal{D}(A^\partial)) = C(E^\partial)$ for all $\lambda > 0$ (see, for example, Theorem 1.2.6 of [9]), an adaptation of the proof of Proposition 4.9.18 of [9] establishes that $\tilde{Y}_t^\infty = \mu Q_t^\partial$, $t \geq 0$.

We have established that the distributions of the processes $(\tilde{Y}^n)_{n=1}^\infty$ are relatively compact and that the limit of any convergent subsequence of these distributions is necessarily the point mass at the function $(\mu Q_t^\partial)_{t\geq 0}$. Therefore, $(\tilde{Y}^n)_{n=1}^\infty$ converges in distribution to $(\mu Q_t^\partial)_{t\geq 0}$. To complete the proof, we recall that if a sequence of random variables defined on the same probability space converges in distribution to a constant, then the sequence also converges in probability to that constant. □

We can now complete the proof of Theorem 1.2. From Lemma 4.1 it is possible, by passing to a subsequence and relabelling if necessary, to assume that the sequence $(\tilde{Y}^n)_{n=1}^\infty$ converges almost surely in the Skorohod topology to $(\mu Q_t^\partial)_{t\geq 0}$. As the limit process is continuous,

(4.3) $$\tilde{Y}_t^n \to \mu Q_t^\partial \text{ almost surely for all } t \geq 0.$$

By the semigroup property of $(P_t)_{t\geq 0}$, it suffices to prove that $\Sigma(\mu Q_t) \leq \Sigma(\mu)$ for all $t > 0$. It is clear from Lemma 2.2 and (4.3) that

$$\Sigma(\mu Q_t) \leq \liminf_{n\to\infty} \Sigma(\tilde{Y}_t^n), \text{ a.s.}$$

for each $t \geq 0$. On the other hand, it is immediate by construction and (4.1) that

$$\Sigma(\tilde{Y}_0^1) \leq \Sigma(\tilde{Y}_0^2) \leq \ldots \leq \Sigma(\mu), \text{ a.s.}$$

Therefore, by part (iv) of Lemma 3.1,

$$\Sigma(\mu Q_t) \leq \liminf_{n\to\infty} \Sigma(\tilde{Y}_t^n) \leq \liminf_{n\to\infty} \Sigma(\tilde{Y}_0^n) \leq \Sigma(\mu).$$

## 5. A COUNTEREXAMPLE

It is not possible to dispense with the assumption in Theorem 1.2 that $W$ is a diffusion, as the following counterexample shows. Let $W$ be the space-time process with $I = \mathbb{R}$ such that for any $(s, x) \in \mathbb{R}_+ \times \mathbb{R}$ we have $\mathbb{P}^{s,x}\{\sigma > t\} = \exp(-t)$ for $\sigma = \inf\{t \geq 0 : X_t \neq X_0\}$ and $\mathbb{P}^{s,x}\{X_\sigma = X_0 \pm 1\} = \frac{1}{2}$ (in particular, $\mathbb{P}^{s,x}\{\zeta = \infty\} = 1$). Note that the semigroup of $W$ is indeed Feller. Take $\mu = \delta_0 - \delta_{1/2}$. It is clear for $t > 0$ that $(\mu Q_t)^+ = \mu^+ Q_t$ is a measure with support the whole of $\mathbb{Z}$ and $(\mu Q_t)^- = \mu^- Q_t$ is a measure with support the whole of $\mathbb{Z} + 1/2$. Thus $1 = \Sigma(\mu) \not\geq \Sigma(\mu Q_t) = \infty$ for $t > 0$.

**Acknowledgement.** The authors thank Steve Altschuler and Lani Wu for helpful discussions. The research was conducted while both authors were visiting the Mathematical Sciences Research Institute, Berkeley, during the Special Year on Stochastic Analysis, 1997–1998.

STEVEN N. EVANS, DEPARTMENT OF STATISTICS #3860, UNIVERSITY OF CALIFORNIA AT BERKELEY, 367 EVANS HALL, BERKELEY, CA 94720-3860
*E-mail address*: evans@stat.berkeley.edu

RUTH J. WILLIAMS, DEPARTMENT OF MATHEMATICS, UNIVERSITY OF CALIFORNIA AT SAN DIEGO, 9500 GILMAN DRIVE, LA JOLLA, CA 92093-0112, U.S.A.
*E-mail address*: williams@math.ucsd.edu
*Current address*: Mathematical Sciences Research Institute, 1000 Centennial Drive, Berkeley, CA 94720-5070